\documentclass[11pt]{article}
\newcommand{\documentdate}{6 XII 2021}

\usepackage{a4wide,latexsym,graphicx,epstopdf,amsmath,varioref,xcolor,lscape}

\title{{\sf OPM}, a collection of Optimization Problems in Matlab}

\author{
  S. Gratton%
  \thanks{Universit\'e de Toulouse, INP, IRIT, Toulouse, France.
    Email: serge.gratton@enseeiht.fr Work partially supported by 3IA Artificial and
     Natural Intelligence Toulouse Institute, French "Investing for the Future
     - PIA3" program under the Grant agreement ANR-19-PI3A-0004"}
  ~and Ph. L. Toint%
  \thanks{NAXYS, University of Namur, Namur, Belgium.
    Email: philippe.toint@unamur.be}
}

\newcommand{\beqn}[1]{\begin{equation}\label{#1}}
\newcommand{\eeqn}{\end{equation}}
\newcommand{\req}[1]{(\ref{#1})}

\newtheorem{theorem}{Theorem}[section]
\newtheorem{lemma}[theorem]{Lemma}
\newtheorem{corollary}[theorem]{Corollary}
\newcommand{\numsection}[1]{\section{#1}\setcounter{equation}{0}}

\newcounter{algo}[section]

\newcommand{\ii}[1]{\{ 1, \ldots, #1 \}}
 
\newcommand{\calE}{{\cal E}}

\newcommand{\smallRe}{\hbox{\footnotesize I\hskip -2pt R}}

\newcommand{\comment}[1]{}
\newcommand{\OPM}{{\sf OPM}}
\newcommand{\mystack}[2]{_{\stackrel{\scriptstyle #1}{\scriptstyle #2}}}

\date{\documentdate}

\begin{document}

\maketitle

\begin{abstract}
\OPM\ is a small collection of {\sf CUTEst} unconstrained and
bound-constrained nonlinear optimization problems, which can be used in Matlab for
testing optimization algorithms directly (i.e.\ without installing additional
software).

\noindent
\textcolor{red}{ == This collection is now superseded by the S2MPJ
  collection (arXiv:2407.07812) ==}
\end{abstract}


\numsection{Introduction}

The purpose of this short report is to introduce \OPM, a small collection of
unconstrained and bound-constrained nonlinear optimization problems,
suitable for testing optimization algorithms written in Matlab\footnote{%
Matlab\textregistered\ is a trade mark of MathWorks, Inc.}. While other collections of this type
are available (we obviously think of the extensive and widely used {\sf CUTEst}
\cite{GoulOrbaToin15b} collection but also of
\cite{MoreGarbHill81,Toin83a,Buck89,BondBortMore99}), none of these is
directly accessible in Matlab, and their interface with Matlab code is often
system dependent.  The collection presented here has no ambition of being
as extensive or as functional as {\sf CUTEst}, which we believe remains an
authoritative source of test problems today, but it may help Matlab
optimization developpers by providing a quick installation-free access to a
small fraction of the {\sf CUTEst} problems and functionalities.

\numsection{The collection}

The \OPM\ collection consists of
\begin{enumerate}
\item a set of Matlab executable files (*.m), each file corresponding to a
  specific problem of the {\sf CUTEst} collection (for instance, the {\tt
    rosenbr.m} file corresponds to ROSENBR.SIF, the famous Rosenbrock banana
  function); 
\item the {\tt opm\_eval\_cpsf.m} file, which provides the mechanism for using
  the other files;
\end{enumerate}
These files are available on \fbox{\tt https://github.com/gratton7/OPM}.  Each problem files starts with a brief
description of the problem and a reference to its source (both extracted from
the corresponding {\sf CUTEest} SIF file).

At the time of writing, the collection contains 134 available problems in
continuous variables that unconstrained or bound constrained (including problems with fixed
variables), i.e.,
\[
\min_{\mystack{x\in\smallRe^n}{x_{\rm low}\leq x \leq x_{\rm up}}} f(x).
\]
It includes the problems given in \cite{MoreGarbHill81} and
\cite{Buck89} as well as selected problems from \cite{Toin83a} and
\cite{BondBortMore99}, all of which can also be found in {\sf CUTEst} \cite{GoulOrbaToin15b}.

\numsection{A simple access to the \OPM\ problem files}

In its simplest form the use of the \OPM\ problems consists of two stages, which
we illustrate by using the {\tt rosenbr} problem.

\begin{enumerate}
\item The first stage is to specify the problem dimension (when possible),
  obtain the standard starting point and, if relevant, the lower and upper
  bounds on the problems' variables.  This achieved by simply inserting the
  following instruction is the user's Matlab program
  \begin{quote}\tt
  x0 = rosenbr( 'setup',10 );
  \end{quote}
  In this case, the problem dimension (that is the number of variables) is set
  to 10.  The second input argument is optional: the simpler call
  \begin{quote}\tt
  x0 = rosenbr( 'setup' );
  \end{quote}
  assigns the default dimension to the problem (2 for {\tt rosenbr}).  Setting
  the problem dimension is possible for the {\tt rosenbr} problem, but not
  necesarily for every problem. When possible, there may be additional
  problem-dependent constraints on how dimension can be chosen\footnote{Check
  the problem m-file to verify the acceptable choices.}.
\item  The value of the objective function, gradient and Hessian at a given
  vector {\tt x} for the Rosenbrock problem are then computed, in the course
  of the user's optimization code, by
  \begin{quote}\tt
  fx = rosenbr( 'objf', x );
  \end{quote}
  if only the objective-function value at {\tt x} is requested, or
  \begin{quote}\tt
  [ fx, gx ] = rosenbr( 'objf', x );
  \end{quote}
  if the objective-function value and its gradient at {\tt x} are requested,
  or finally
  \begin{quote}\tt
  [ fx, gx, Hx ] = rosenbr( 'objf', x );
  \end{quote}
  if value, gradient and Hessian at {\tt x} are requested. If the above, {\tt
    fx} is a scalar, {\tt x} and {\tt gx} are column vectors of size $n$, say,
  and {\tt Hx} is an $n \times n$ square matrix.
\end{enumerate}

There are variations on this simple calling sequences. In particular, bounds
on the variables can also be retrieved, as we explain in the next section.

\numsection{An access to problems exploiting their
  coordinate-partially-separable structure}

All OPM test problems are, just as a
substantial fraction of real problems from applications,
``coordinate-partially-separable'' (CPS). A CPS problem is a problem whose
objective function may be written in the form
\beqn{cpsf}
f(x) = \sum_{i=1}^{n_e} f_i( x_{\calE_i}),
\eeqn
where the sets $\calE_i$ are subsets of $\ii{n}$ and
$x_{\calE_i}$ is the subvector of $x$ indexed by $\calE_i$.  In other words,
$f$ is the sum of $n_e$ \emph{element functions} $f_i$, each one of them
only depending on a subset of the variables. For large problems, it very often
the case that $|\calE_i|$, the size of $\calE_i$, is independent of the
problem dimension and $\max_i|\calE_i| \ll n$, as it can be shown
\cite{GrieToin82a} that every sufficiently smooth function whose Hessian
matrix is sparse if a CPS function. (Indeed, the $|\calE_i|$ are the
dimensions of the Hessian's dense principal submatrices.) In \req{cpsf}, the
$f_i$ are called the \emph{element functions}, the $\calE_i$ the
\emph{element domains} and \req{cpsf} the \emph{CPS decomposition} of the
function $f$.

A quick look at the {\tt opm\_eval\_cpsf.m} file reveals that this Matlab function
merely builds the sum \req{cpsf}. But \OPM\ also provides the means to access
the element functions $f_i$ (and their gradients and Hessians)
individually. For example, the following sequence of calls applies to the {\tt
  lminsurf} problem: one starts by setting up the problem and retrieving its
structure in a Matlab {\tt struct} called {\tt cpsstr}\footnote{For CPS
  structure.}. This is done in the following calling sequence:
\begin{quote}\tt
  x0     = lminsurf( 'setup' );\\
  cpsstr = lminsurf( 'cpsstr', length(x0) );
\end{quote}
In particular, its element domains are given in the cell {\tt cpsstr.eldom}
whose $i$-th entry is a vector containing the indices of $\calE_i$.  Other
structure parameters are given in the cell {\tt cpsstr.param}, when relevant.
Suppose now that, in the course of the user's program, the value
of the {\tt i}-th element function (for {\tt i} $\in \ii{n_e}$), its gradient
and Hessian are needed at the vector of variables {\tt x}.  These values are
simply obtained by the call
\begin{quote}\tt
  [ fix, gix, Hix ] = ...\\
   \hspace*{1cm} lminsurf( 'elobjf', i, x( cpsstr.eldom\{i\} ), cpsstr.param\{:\} );
\end{quote}
Note that {\tt gix}
is a column vector of length $|\calE_i| = ${\tt length(cpsstr.eldom\{i\})} giving the first
derivatives of $f_i$ with respect to the variables whose index is in
$\calE_i$.  Similarly, {\tt Hix} is a symmetric square matrix of dimension
$|\calE_i|\times|\calE_i|$ containing the second derivatives of $f_i$ with
respect to these variables.  If the $i$-th element function involves
$p$ numerical parameters {\tt par1} to {\tt parp}, a call of the form
\begin{quote}\tt
  [ fix, gix, Hix ] = ...\\
  \hspace*{1cm} lminsurf( 'elobjf', i, x( cpsstr.eldom\{i\} ),\\
  \hspace*{2cm} cpsstr.param\{:\}, par1,\ldots, parp  );
\end{quote}
is allowed.

The use of the CPS structure may be numerically extremely advantageous for large
problems where $\max_i|\calE_i| \ll n$, as it gives direct access to the sparsity
structure of the objective function's Hessian. It is also possible to use the
element gradients to construct 'elementwise quasi-Newton'
methods\footnote{Such ``partioned updating'' technique were proposed in
\cite{GrieToin82b} and motivated the introduction of partially separable
functions.}. In model-based derivative-free algorithms, element-wise models
can be used with considerable success \cite{ColsToin05,PorcToin21}.

\numsection{An alternative calling sequence}

The reader might wonder if obtaining the problem structure by calling
the problem with first argument '{\tt cpsstr}', as explained in
the previous section, might be useful even if access to individual element
functions is not desired.  And (you guessed it), the answer is positive.
Instead to the simple calls
\begin{quote}\tt
  x0 = rosenbr( 'setup',10 );\\
  ...\\
  {}[ fx, gx, Hx ] = rosenbr( 'objf', x );\\
\end{quote}
may be replaced by the barely more complex
\begin{quote}\tt
  x0     = rosenbr( 'setup', 10 );\\
  cpsstr = rosenbr( 'cpsstr', 10 );\\
  ...\\
  {}[ fx, gx, Hx ] = rosenbr( 'objf', x, cpsstr );\\
\end{quote}
In effect, this last calling sequence computes the problem's structure (in
{\tt cpsstr}) just once, and then passes it to every subsequent call for
calculating the objective function's value (and derivatives). By contrast, the
first calling sequence recomputes the problem structure each time a 
call for calculating the objective function's value (and derivatives) is made.
Depending on problems, this may slow down evaluations.

\numsection{A formal description of the calling sequence}

The simple calling sequence discussed above is a special case of what is
possible. We now give the complete list of arguments for the calls to a
problems file.  Suppose we consider the problem {\tt problem} ({\tt rosenbr}
above).  Then the call to {\tt problem} always has the form
\begin{quote}\tt
  varargout = problem( action, varargin );
\end{quote}
where
\begin{description}
\item[\tt action:] is a string whose possible values are
  \begin{description}
  \item[\tt 'setup'] when the call is to set up the problem,
  \item[\tt 'eldom'] when the lement domains of the problem are requested,
  \item[\tt 'objf'] when the call requests the value of the problem's objective
    function (and possibly those of its gradient and Hessian) at a given vector of
    variables,
  \item[\tt 'elobjf'] when the call requests the value of a specific element
    function of the problem (and possibly those of its gradient and Hessian)
    at a given vector of variables,
  \item[\tt 'consf'] when the call requests values of the
    constraint's function (and possibly their Jacobian and Hessian) at a given
    vector of variables.
  \end{description}
\item[\tt varargin:] is an optional list of problems parameters.
\end{description}

\vspace*{5mm}
\noindent
$\bullet$ For a \emph{setup call}, the complete calling sequence is given by
\begin{quote}\tt
[ x0, fstar, xtype, xlower, xupper, clower, cupper, class ] ...\\
\hspace*{1cm}  = problem( 'setup', varargin )
\end{quote}
where, on input,
\begin{description}
\item[\tt varargin:] is an optional list of numerical parameters for the problem. When specified,
  the first element of the list must be a positive integer and is interpreted
  as the problem's dimension.
\end{description}
For a problem of dimension $n$, the meaning of the ouput arguments is as follows.
\begin{description}
\item[\tt x0:] is a column vector of dimension $n$ giving the standard problem's starting point ;
\item[\tt fstar:] is either a numerical value giving the value of the
  objective function at a known minimizer, or a string giving some information
  on this value, when possible;
\item[\tt xtype:] is a string of length $n$ whose $i$-th character indicates
  the type of the $i$-th variable, and can take the values
  \begin{description}
  \item[\tt 'c':] for a continuous variable,
  \item[\tt 'i':] for an integer variable,
  \item[\tt 's':] for a categorical variable;
  \end{description}
  A empty string is equivalent to a string of all {\tt 'c'}, indicating that
  all variables are continuous. (Note that only {\tt 'c'} variables are used
  in the current problem set.) 
\item[\tt xlower:] is a column vector of size $n$, giving the lower bounds on the
  variable.  Notice that {\tt -Inf} is an acceptable lower bound value. An
  empty vector is equivalent to a vector whose every component is equal to
  {\tt -Inf}, indicating that none of the variables has a finite lower bound.
\item[\tt xupper:] is a column vector of size $n$, giving the upper bounds on the
  variable.  Notice that {\tt +Inf} is an acceptable upper bound value. An
  empty vector is equivalent to a vector whose every component is equal to
  {\tt +Inf}, indicating that none of the variables has a finite upper bound.
\item[\tt clower:] is a column vector of size equal to the number of constraints,
  giving the lower bounds on the value of the constraint functions.
  Notice that {\tt -Inf} is an acceptable lower bound value. An
  empty vector is equivalent to a vector whose every component is equal to
  {\tt -Inf}, indicating that all constraints are of the ``less than'' form.
\item[\tt cupper:] is a column vector of size equal to the number of constraints,
  giving the upper bounds on the value of the constraint function.
  Notice that {\tt +Inf} is an acceptable upper bound value. An
  empty vector is equivalent to a vector whose every component is equal to
  {\tt +Inf}, indicating that all constraints are of the ``greater than'' form.
\item[\tt class:] is a string indicating the class of the problem, according to
  the {\sf CUTEst} classification scheme (see \cite{GoulOrbaToin15b}).
\end{description}
Note that, as usual in Matlab, incomplete lists of output arguments are
allowed.

\vspace*{5mm}
\noindent
$\bullet$ For a \emph{CPS structure call}, the calling sequence reduces to
\begin{quote}\tt
  cpsstr = problem( 'cpsstr', n );
\end{quote}
where, on input,
\begin{description}
\item[\tt n:] is the problem's dimension,
\end{description}
while, on output
\begin{description}
\item[\tt cpsstr:] is a Matlab struct with fields
  \begin{description}
  \item[\tt name:] a string giving the problem's name,
  \item[\tt eldom:] a cell of length equal to the number $n_e$ of element
  functions in \req{cpsf} and whose $i$-th entry is a vector containing the
  indices of $\calE_i$.,
  \item[\tt param:] a cell containing additional problem parameters.  For
    instance, this might contain a vector of measurements for a nonlinear
    least-squares fitting problem, or, more simply, the problem's dimension.
    This cell may be empty if no parameter is required for evaluating either
    the full objective function or its CPS element functions.
  \end{description}
\end{description}
      
\vspace*{5mm}
\noindent
$\bullet$ For an \emph{full objective value call} requiring objective-function (and possibly derivatives) values, the
calling sequence is
\begin{quote}\tt
  [ fx, gx, Hx ] = problem( 'objf', x, varargin );
\end{quote}
or
\begin{quote}\tt
  [ fx, gx, Hx ] = problem( 'objf', x, cpsstr, varargin );
\end{quote}
where, on input,
\begin{description}
\item[\tt x:] is a column vector of dimension $n$, specifying the point at
  which the values must be computed,
\item[\tt cpsstr:] is the {\tt struct} describing the problem's CPS
  structure, as given on output of a (previous) CPS structure call,
\item[\tt varargin:] is an optional list of numerical parameters for the problem,
\end{description}
and, on output,
\begin{description}
\item[\tt fx:] is the returned computed objective-function value at {\tt x},
\item[\tt gx:] is a column vector of dimension $n$ containing the returned
  computed objective-function gradient at {\tt x}, 
\item[\tt Hx:] is a square symmetric matrix of dimension $n \times n$
  containing the returned computed objective-function Hessian at {\tt x}.
\end{description}
Note again that, as usual in Matlab, incomplete list of output arguments are
allowed. Moreover, an output which is not requested is not calculated.  Thus a call
\begin{quote}\tt
  fx = problem( 'objf', x, varargin );
\end{quote}
only computes the objective's value {\tt fx} (the gradient and Hessian are not
computed), and, similarly the call
\begin{quote}\tt
  [ fx, gx ] = problem( 'objf', x, varargin );
\end{quote}
does not computes the Hessian. Both calls however recompute the problem's
structure, which is not the case for the calls
\begin{quote}\tt
  fx = problem( 'objf', x, cpsstr, varargin );
\end{quote}
and
\begin{quote}\tt
  fx = problem( 'objf', x, cpsstr, varargin );
\end{quote}

\vspace*{5mm}
\noindent
$\bullet$ For an \emph{element objective value call} requiring
the value (and possibly derivatives values) of a specific element function
of the CPS decomposition \req{cpsf}, the calling sequence is
\begin{quote}\tt
  [ fix, gix, Hix ] = ...\\
  \hspace*{1cm} problem( 'elobjf', i, xi, cpsstr.param\{:\}, varargin );
\end{quote}
where, on input,
\begin{description}
\item[\tt i:] is the index (in $\{1, \ldots, n_e\}$) of the considered
  element function;
\item[\tt xi:] is a column vector of dimension $|\calE_i|$, specifying the
  values of the variables occuring in the {\tt i}-th element function (as
  specified in {\tt cpsstr.eldom\{i\}}) for which the outputs of this
  element function must be calculated,
\item[\tt cpsstr:] is the {\tt struct} describing the problem's CPS
  structure, as given on output of a (previous) CPS structure call (see above),
\item[\tt varargin:] is an optional list of numerical parameters for the problem,
\end{description}
and, on output,
\begin{description}
\item[\tt fix:] is the returned computed element objective-function value
  $f_i$ at {\tt xi},
\item[\tt gix:] is a column vector of dimension $|\calE_i|$ containing the returned
  computed element objective-function gradient $\nabla_x^1f_i$ at {\tt xi}, with respect to the
  variables specified in  {\tt cpsstr.eldom\{i\}},
\item[\tt Hix:] is a square symmetric matrix of dimension $|\calE_i| \times |\calE_i|$
  containing the returned computed element objective-function Hessian
  $\nabla_x^2f_i$ at {\tt xi}, with respect to the variables specified in
  {\tt cpsstr.eldom\{i\}}.
\end{description}
Note once more that an incomplete list of output arguments is
allowed.

\numsection{Some provisions for future extensions}

The reader has undoubtly noticed that some possibilities have not been fully
described here, notably variable types other than continuous and
evaluation of explicit constraint functions.  There are included at this stage
merely to provide room for future development, should this prove to be useful.
  
\numsection{Conclusion}

We have briefly outlined the \OPM\ optimization test problem collection and
provided guidance on how to use it.  The authors of course welcome
contributions from users in the form of additional (duly verified) test
problems\footnote{If you have one, please contact one of the authors.}.  Other
future developments are contingent on users' demand (and developpers' time).

{\footnotesize

}

\newpage
\appendix

\vspace*{5mm}
\noindent
\textbf{\Large Problem descriptions}

\vspace*{5mm}
\noindent
In the following pages, we give a summary description of each problem, where
the meaning of the column's headers are as follows:
\begin{center}
\begin{tabular}{ll}
{\tt fstar:} & the value of the objective function at a minimizer (if known);\\
{\tt n:}& the problem's default dimension;\\
{\tt mel:}& the maximum element-domain size ($\max_i|\calE_i|$),\\
{\tt nc} & the problem's number of continuous variables;\\
{\tt ni:} & the problem's number of integer variables;\\
{\tt nfree:} & the problem's number of free variables;\\
{\tt nlow:} & the problem's number of variables that are bounded below;\\
{\tt nupp:} & the problem's number of variables that are bounded above;\\
{\tt nfix:} & the problem's number of variables that are fixed;\\
{\tt m:} & the problem's number of explicit constraints;\\
{\tt mi:} & the problem's number of explicit inequality constraints;\\
{\tt me:} & the problem's number of explicit equality constraints;\\
{\tt objtype:} & some information on the problem's objective function;\\
{\tt constype:} & some information on the problem's constraints;\\
{\tt classif:} & the problem's {\sf CUTEst} classification string.
\end{tabular}
\end{center}

\newpage
\begin{landscape}
{\footnotesize
\begin{verbatim}
**************************************************************************************************************************
**************************************************************************************************************************
********************************                                           ***********************************************
********************************     OPM library content (17-Nov-2021)     ***********************************************
********************************                                           ***********************************************
**************************************************************************************************************************
**************************************************************************************************************************
 
Problem name        fstar         n    mel  nc   ni nfree  nlow  nupp nboth  nfix  m   mi   me objtype constype  classif
 
argauss         +1.12793277e-08    3    3    3    0    3     0     0    0     0    0    0    0  nonlin     -   SUR2-AY-3-0     
arglina         +1.00000000e+01   10   10   10    0   10     0     0    0     0    0    0    0  nonlin     -   SUR2-AN-V-0     
arglinb         +4.87179487e+00   10   10   10    0   10     0     0    0     0    0    0    0  nonlin     -   SUR2-AN-V-0     
arglinc         +6.13513514e+00   10   10   10    0   10     0     0    0     0    0    0    0  nonlin     -   SUR2-AN-V-0     
argtrig         +0.00000000e+00   10   10   10    0   10     0     0    0     0    0    0    0  nonlin     -   SUR2-AN-V-0     
arwhead         +0.00000000e+00   10    2   10    0   10     0     0    0     0    0    0    0  nonlin     -   OUR2-AN-V-0     
bard            +8.21500000e-03    3    3    3    0    3     0     0    0     0    0    0    0  nonlin     -   SUR2-AY-3-0     
bdarwhd         +0.00000000e+00   10    3   10    0   10     0     0    0     0    0    0    0  nonlin     -   OUR2-AY-V-0     
beale           +0.00000000e+00    2    2    2    0    2     0     0    0     0    0    0    0  nonlin     -   SUR2-AN-2-0     
biggs5          +0.00000000e+00    6    6    6    0    5     0     0    0     1    0    0    0  nonlin     -   SXR2-AN-6-0     
biggs6          +0.00000000e+00    6    6    6    0    6     0     0    0     0    0    0    0  nonlin     -   SUR2-AN-6-0     
booth           +0.00000000e+00    2    2    2    0    2     0     0    0     0    0    0    0  nonlin     -   QUR2-AN-2-0     
box2            +0.00000000e+00    3    3    3    0    2     0     0    0     1    0    0    0  nonlin     -   OXR2-AN-2-0     
box3            +0.00000000e+00    3    3    3    0    3     0     0    0     0    0    0    0  nonlin     -   OUR2-AN-3-0     
bratu1d             unknown       30    2   30    0   28     0     0    0     2    0    0    0  nonlin     -   SUR2-AY-V-0     
brkmcc               -Inf          2    2    2    0    2     0     0    0     0    0    0    0  nonlin     -   OUR2-AN-2-0     
brownal         +0.00000000e+00   10   10   10    0   10     0     0    0     0    0    0    0  nonlin     -   SUR2-AN-V-0     
brownbs         +0.00000000e+00    2    2    2    0    2     0     0    0     0    0    0    0  nonlin     -   SUR2-AN-2-0     
brownden        +8.58222000e+04    4    4    4    0    4     0     0    0     0    0    0    0  nonlin     -   SUR2-AN-4-0     
broyden3d       +0.00000000e+00   10    3   10    0    8     0     0    0     2    0    0    0  nonlin     -   SXR2-AY-V-0     
broydenbd       +0.00000000e+00   10    7   10    0   10     0     0    0     0    0    0    0  nonlin     -   SUR2-AY-V-0     
chandheu        +0.00000000e+00   57   57   57    0   57     0     0    0     0    0    0    0  nonlin     -   SUR2-RN-V-0     
chebyqad        +0.00000000e+00   10   10   10    0   10     0     0    0     0    0    0    0  nonlin     -   SUR2-AN-V-0     
cliff           +1.99786610e-01    2    2    2    0    2     0     0    0     0    0    0    0  nonlin     -   OUR2-AN-2-0     
clplatea            unknown       16    2   16    0   12     0     0    0     4    0    0    0  nonlin     -   OXR2-MN-V-0     
clplateb            unknown       16    4   16    0   12     0     0    0     4    0    0    0  nonlin     -   OXR2-MN-V-0     
clustr          +0.00000000e+00    2    2    2    0    2     0     0    0     0    0    0    0  nonlin     -   SUR2-AN-2-0     
cosine          +0.00000000e+00   10    2   10    0   10     0     0    0     0    0    0    0  nonlin     -   QUR2-AY-V-0     
crglvy          +0.00000000e+00   10    2   10    0   10     0     0    0     0    0    0    0  nonlin     -   SUR2-AY-V-0     
cube            +0.00000000e+00   10    2   10    0   10     0     0    0     0    0    0    0  nonlin     -   SUR2-AY-V-0     

**************************************************************************************************************************
**************************************************************************************************************************
********************************                                           ***********************************************
********************************     OPM library content (17-Nov-2021)     ***********************************************
********************************                                           ***********************************************
**************************************************************************************************************************
**************************************************************************************************************************
 
Problem name        fstar         n    mel  nc   ni nfree  nlow  nupp nboth  nfix  m   mi   me objtype constype  classif
 
curly10             unknown       30   11   30    0   30     0     0    0     0    0    0    0  nonlin     -   OUR2-AN-V-0     
curly20             unknown       30   21   30    0   30     0     0    0     0    0    0    0  nonlin     -   OUR2-AN-V-0     
curly30             unknown       40   31   40    0   40     0     0    0     0    0    0    0  nonlin     -   OUR2-AN-V-0     
deconvu              -Inf         51   22   51    0   51     0     0    0     0    0    0    0  nonlin     -   SUR2-MN-51-0    
dixmaana        +1.00000000e+00   12    4   12    0   12     0     0    0     0    0    0    0  nonlin     -   OUR2-AY-V-0     
dixmaanb        +1.00000000e+00   12    4   12    0   12     0     0    0     0    0    0    0  nonlin     -   OUR2-AY-V-0     
dixmaanc        +1.00000000e+00   12    4   12    0   12     0     0    0     0    0    0    0  nonlin     -   OUR2-AY-V-0     
dixmaand        +1.00000000e+00   12    4   12    0   12     0     0    0     0    0    0    0  nonlin     -   OUR2-AY-V-0     
dixmaane        +1.00000000e+00   12    4   12    0   12     0     0    0     0    0    0    0  nonlin     -   OUR2-AY-V-0     
dixmaanf        +1.00000000e+00   12    4   12    0   12     0     0    0     0    0    0    0  nonlin     -   OUR2-AY-V-0     
dixmaang        +1.00000000e+00   12    4   12    0   12     0     0    0     0    0    0    0  nonlin     -   OUR2-AY-V-0     
dixmaanh        +1.00000000e+00   12    4   12    0   12     0     0    0     0    0    0    0  nonlin     -   OUR2-AY-V-0     
dixmaani        +1.00000000e+00   12    4   12    0   12     0     0    0     0    0    0    0  nonlin     -   OUR2-AY-V-0     
dixmaanj        +1.00000000e+00   12    4   12    0   12     0     0    0     0    0    0    0  nonlin     -   OUR2-AY-V-0     
dixmaank        +1.00000000e+00   12    4   12    0   12     0     0    0     0    0    0    0  nonlin     -   OUR2-AY-V-0     
dixmaanl        +1.00000000e+00   12    4   12    0   12     0     0    0     0    0    0    0  nonlin     -   OUR2-AY-V-0     
dixon           +0.00000000e+00   10    2   10    0   10     0     0    0     0    0    0    0  nonlin     -   QUR2-AY-V-0     
dqrtic          +0.00000000e+00   10    1   10    0   10     0     0    0     0    0    0    0  nonlin     -   QUR2-AY-V-0     
edensch             unknown       10    2   10    0   10     0     0    0     0    0    0    0  nonlin     -   OUR2-AY-V-0     
eg2                 unknown       10    1   10    0   10     0     0    0     0    0    0    0  nonlin     -   OUR2-AY-V-0     
eg2s                unknown       10    3   10    0   10     0     0    0     0    0    0    0  nonlin     -   OUR2-AY-V-0     
eigenals        +0.00000000e+00  110   30  110    0  110     0     0    0     0    0    0    0  nonlin     -   SUR2-AN-V-0     
eigenbls        +0.00000000e+00  110   30  110    0  110     0     0    0     0    0    0    0  nonlin     -   SUR2-AN-V-0     
eigencls        +0.00000000e+00  110   30  110    0  110     0     0    0     0    0    0    0  nonlin     -   SUR2-AN-V-0     
engval1         +0.00000000e+00   10    2   10    0   10     0     0    0     0    0    0    0  nonlin     -   OUR2-AY-V-0     
engval2         +0.00000000e+00    3    3    3    0    3     0     0    0     0    0    0    0  nonlin     -   SUR2-AY-V-0     
expfit          +0.00000000e+00    2    2    2    0    2     0     0    0     0    0    0    0  nonlin     -   SUR2-AN-2-0     
extrosnb        +0.00000000e+00   10    2   10    0   10     0     0    0     0    0    0    0  nonlin     -   SUR2-AY-V-0     
fminsurf        +1.00000000e+00   16   16   16    0   16     0     0    0     0    0    0    0  nonlin     -   OUR2-MY-V-0     

**************************************************************************************************************************
**************************************************************************************************************************
********************************                                           ***********************************************
********************************     OPM library content (17-Nov-2021)     ***********************************************
********************************                                           ***********************************************
**************************************************************************************************************************
**************************************************************************************************************************
 
Problem name        fstar         n    mel  nc   ni nfree  nlow  nupp nboth  nfix  m   mi   me objtype constype  classif
 
freuroth        +0.00000000e+00   10    2   10    0   10     0     0    0     0    0    0    0  nonlin     -   OUR2-AY-V-0     
genhumps        +0.00000000e+00   10    2   10    0   10     0     0    0     0    0    0    0  nonlin     -   SUR2-AY-V-0     
gottfr          +0.00000000e+00    2    2    2    0    2     0     0    0     0    0    0    0  nonlin     -   SUR2-AN-2-0     
gulf            +0.00000000e+00    3    3    3    0    3     0     0    0     0    0    0    0  nonlin     -   SUR2-AN-3-0     
hairy           +2.00000000e+01    2    2    2    0    2     0     0    0     0    0    0    0  nonlin     -   OUR2-AN-2-0     
heart6ls        +2.00000000e+00    6    6    6    0    6     0     0    0     0    0    0    0  nonlin     -   SUR2-AY-6-0     
heart8ls        +0.00000000e+00    8    8    8    0    8     0     0    0     0    0    0    0  nonlin     -   SUR2-MN-8-0     
helix           +0.00000000e+00   10    3   10    0   10     0     0    0     0    0    0    0  nonlin     -   SUR2-AY-V-0     
hilbert         +0.00000000e+00   10   10   10    0   10     0     0    0     0    0    0    0  nonlin     -   QUR2-AN-V-0     
himln3          -1.00000000e+00    2    2    2    0    2     0     0    0     0    0    0    0  nonlin     -   OUR2-AN-2-0     
himm25          +0.00000000e+00    2    2    2    0    2     0     0    0     0    0    0    0  nonlin     -   SUR2-AN-2-0     
himm27          +0.00000000e+00    2    2    2    0    2     0     0    0     0    0    0    0  nonlin     -   OUR2-AN-2-0     
himm28          +0.00000000e+00    2    2    2    0    2     0     0    0     0    0    0    0  nonlin     -   SUR2-AN-2-0     
himm29          +0.00000000e+00    2    2    2    0    2     0     0    0     0    0    0    0  nonlin     -   SUR2-AN-2-0     
himm30          +0.00000000e+00    3    3    3    0    3     0     0    0     0    0    0    0  nonlin     -   SUR2-AN-3-0     
himm32          +3.18572000e+02    4    4    4    0    4     0     0    0     0    0    0    0  nonlin     -   SUR2-AN-4-0     
himm33          +0.00000000e+00    2    2    2    0    2     0     0    0     0    0    0    0  nonlin     -   OUR2-AN-2-0     
hydc20ls        +0.00000000e+00   99   14   99    0   99     0     0    0     0    0    0    0  nonlin     -   SUR2-AN-99-0    
hypcir          +1.23259516e-01    2    2    2    0    2     0     0    0     0    0    0    0  nonlin     -   SUR2-AN-2-0     
indef               unknown       10    3   10    0   10     0     0    0     0    0    0    0  nonlin     -   OUR2-AY-V-0     
integreq        +0.00000000e+00   10   10   10    0   10     0     0    0     0    0    0    0  nonlin     -   SUR2-AN-V-0     
jensmp          +1.24362000e+02    2    2    2    0    2     0     0    0     0    0    0    0  nonlin     -   SUR2-AN-2-0     
kowosb          +3.07505000e-03    4    4    4    0    4     0     0    0     0    0    0    0  nonlin     -   SUR2-AN-4-0     
lminsurf        +9.00000000e+00   16    4   16    0    4     0     0    0    12    0    0    0  nonlin     -   OXR2-MY-V-0     
mancino             unknown       10   10   10    0   10     0     0    0     0    0    0    0  nonlin     -   SUR2-AN-V-0     
mexhat          -1.11715260e+00    2    2    2    0    2     0     0    0     0    0    0    0  nonlin     -   SUR2-AN-2-0     
meyer3          +8.79458000e+01    3    3    3    0    3     0     0    0     0    0    0    0  nonlin     -   SUR2-RN-3-0     
morebv          +0.00000000e+00   12    3   12    0   10     0     0    0     2    0    0    0  nonlin     -   SXR2-AY-V-0     
msqrtals        +0.00000000e+00   16    7   16    0   16     0     0    0     0    0    0    0  nonlin     -   SUR2-AY-V-0     

**************************************************************************************************************************
**************************************************************************************************************************
********************************                                           ***********************************************
********************************     OPM library content (17-Nov-2021)     ***********************************************
********************************                                           ***********************************************
**************************************************************************************************************************
**************************************************************************************************************************
 
Problem name        fstar         n    mel  nc   ni nfree  nlow  nupp nboth  nfix  m   mi   me objtype constype  classif
 
msqrtbls        +0.00000000e+00   16    7   16    0   16     0     0    0     0    0    0    0  nonlin     -   SUR2-AY-V-0     
ncb20               unknown       30   30   30    0   30     0     0    0     0    0    0    0  nonlin     -   OUR2-AN-V-0     
ncb20b              unknown       21   20   21    0   21     0     0    0     0    0    0    0  nonlin     -   OUR2-AN-V-0     
ncb20c              unknown       30   30   30    0   30     0     0    0     0    0    0    0  nonlin     -   OUR2-AN-V-0     
nlminsurf           unknown       16    4   16    0    4     0     0    0    12    0    0    0  nonlin     -   OXR2-MY-V-0     
nondia          +0.00000000e+00   10    2   10    0   10     0     0    0     0    0    0    0  nonlin     -   SUR2-AY-V-0     
nondquar        +0.00000000e+00  100    3  100    0  100     0     0    0     0    0    0    0  nonlin     -   OUR2-AN-V-0     
nzf1            +0.00000000e+00   13    6   13    0   13     0     0    0     0    0    0    0  nonlin     -   OUR2-AY-V-0     
osbornea        +5.46489000e-05    5    5    5    0    5     0     0    0     0    0    0    0  nonlin     -   SUR2-AN-5-0     
osborneb        +4.01377000e-02   11   11   11    0   11     0     0    0     0    0    0    0  nonlin     -   SUR2-AN-11-0    
penalty1            unknown       10   10   10    0   10     0     0    0     0    0    0    0  nonlin     -   SUR2-AY-V-0     
penalty2            unknown       10   10   10    0   10     0     0    0     0    0    0    0  nonlin     -   SUR2-AN-V-0     
penalty3            unknown       10   10   10    0   10     0     0    0     0    0    0    0  nonlin     -   OUR2-AY-V-0     
powellbs        +0.00000000e+00    2    2    2    0    2     0     0    0     0    0    0    0  nonlin     -   OUR2-AN-2-0     
powellsg        +0.00000000e+00    4    2    4    0    4     0     0    0     0    0    0    0  nonlin     -   SUR2-AY-V-0     
powellsq        +0.00000000e+00    2    2    2    0    2     0     0    0     0    0    0    0  nonlin     -   OUR2-AN-2-0     
powr            +0.00000000e+00   10   10   10    0   10     0     0    0     0    0    0    0  nonlin     -   SUR2-AN-V-0     
recipe          +0.00000000e+00    3    3    3    0    3     0     0    0     0    0    0    0  nonlin     -   SUR2-AN-3-0     
rosenbr         +0.00000000e+00    2    2    2    0    2     0     0    0     0    0    0    0  nonlin     -   SUR2-AY-V-0     
s308            +7.73199000e-01    2    2    2    0    2     0     0    0     0    0    0    0  nonlin     -   SUR2-AN-2-0     
schmvett        -3.00000000e+00    3    3    3    0    3     0     0    0     0    0    0    0  nonlin     -   OUR2-AN-3-0     
scosine         +0.00000000e+00   10    2   10    0   10     0     0    0     0    0    0    0  nonlin     -   QUR2-AY-V-0     
scurly10            unknown       30   11   30    0   30     0     0    0     0    0    0    0  bad scal.  -   OUR2-AN-V-0     
scurly20            unknown       30   21   30    0   30     0     0    0     0    0    0    0  bad scal.  -   OUR2-AN-V-0     
scurly30            unknown       30   30   30    0   30     0     0    0     0    0    0    0  bad scal.  -   OUR2-AN-V-0     
sensors             unknown       10   10   10    0   10     0     0    0     0    0    0    0  nonlin     -   OUR2-AN-V-0     
shydc20ls       +0.00000000e+00   99   14   99    0   99     0     0    0     0    0    0    0  nonlin     -   SUR2-AN-99-0    
sisser          +0.00000000e+00    2    2    2    0    2     0     0    0     0    0    0    0  nonlin     -   OUR2-AN-2-0     
spmsqrt         +0.00000000e+00   10    5   10    0   10     0     0    0     0    0    0    0  nonlin     -   SUR2-AY-V-0     

**************************************************************************************************************************
**************************************************************************************************************************
********************************                                           ***********************************************
********************************     OPM library content (17-Nov-2021)     ***********************************************
********************************                                           ***********************************************
**************************************************************************************************************************
**************************************************************************************************************************
 
Problem name        fstar         n    mel  nc   ni nfree  nlow  nupp nboth  nfix  m   mi   me objtype constype  classif
 
tcontact        +2.48015700e+00  100    4  100    0   48     0     0    0    52    0    0    0  nonlin     -   OXR2-MY-V-0     
tlminsurf       +9.00000000e+00   16    4   16    0    4     0     0    0    12    0    0    0  nonlin     -   OXR2-MY-V-0     
tlminsurfx      +9.00000000e+00   16    4   16    0    4     0     0    0    12    0    0    0  nonlin     -   OXR2-MY-V-0     
tnlminsurf          unknown       16    4   16    0    4     0     0    0    12    0    0    0  nonlin     -   OXR2-MY-V-0     
tnlminsurfx         unknown       16    4   16    0    4     0     0    0    12    0    0    0  nonlin     -   OXR2-MY-V-0     
tquartic        +0.00000000e+00   10    1   10    0   10     0     0    0     0    0    0    0  nonlin     -   OUR2-AY-V-0     
tridia          +0.00000000e+00   10    2   10    0   10     0     0    0     0    0    0    0  nonlin     -   QUR2-AY-V-0     
trigger         +0.00000000e+00    7    5    7    0    6     0     0    0     1    0    0    0  nonlin     -   SXR2-AN-7-0     
vardim          +0.00000000e+00   10   10   10    0   10     0     0    0     0    0    0    0  nonlin     -   SUR2-AY-V-0     
vibrbeam            unknown        8    8    8    0    8     0     0    0     0    0    0    0  nonlin     -   SUR2-MN-8-0     
watson              unknown       12   12   12    0   12     0     0    0     0    0    0    0  nonlin     -   SUR2-AN-V-0     
wmsqrtals           unknown       16    7   16    0   16     0     0    0     0    0    0    0  nonlin     -   SUR2-AY-V-0     
wmsqrtbls           unknown       16    7   16    0   16     0     0    0     0    0    0    0  nonlin     -   SUR2-AY-V-0     
woods           +0.00000000e+00   12    2   12    0   12     0     0    0     0    0    0    0  nonlin     -   OUR2-AY-V-0     
yfitu           +0.00000000e+00    3    3    3    0    3     0     0    0     0    0    0    0  nonlin     -   SUR2-MN-3-0     
zangwil2        -1.82000000e+01    2    2    2    0    2     0     0    0     0    0    0    0  nonlin     -   QUR2-AN-2-0     
zangwil3        +0.00000000e+00    3    3    3    0    3     0     0    0     0    0    0    0  nonlin     -   QUR2-AN-3-0     
 
**************************************************************************************************************************
\end{verbatim}
}
\end{landscape}


\begin{thebibliography}{1}

\bibitem{BondBortMore99}
A.~S. Bondarenko, D.~M. Bortz, and J.~J. Mor\'{e}.
\newblock {COPS}: Large-scale nonlinearly constrained optimization problems.
\newblock Technical Report ANL/MCS-TM-237, Mathematics and Computer Science,
  Argonne National Laboratory, Argonne, Illinois, USA, 1999.

\bibitem{Buck89}
A.~G. Buckley.
\newblock Test functions for unconstrained minimization.
\newblock Technical Report CS-3, Computing Science Division, Dalhousie
  University, Dalhousie, Canada, 1989.

\bibitem{ColsToin05}
B.~Colson and Ph.~L. Toint.
\newblock Optimizing partially separable functions without derivatives.
\newblock {\em Optimization Methods and Software}, 20(4-5):493--508, 2005.

\bibitem{GoulOrbaToin15b}
N.~I.~M. Gould, D.~Orban, and Ph.~L. Toint.
\newblock {\sf CUTEst}: a constrained and unconstrained testing environment
  with safe threads for mathematical optimization.
\newblock {\em Computational Optimization and Applications}, 60(3):545--557,
  2015.

\bibitem{GrieToin82a}
A.~Griewank and Ph.~L. Toint.
\newblock On the unconstrained optimization of partially separable functions.
\newblock In M.~J.~D. Powell, editor, {\em Nonlinear Optimization 1981}, pages
  301--312, London, 1982. Academic Press.

\bibitem{GrieToin82b}
A.~Griewank and Ph.~L. Toint.
\newblock Partitioned variable metric updates for large structured optimization
  problems.
\newblock {\em Numerische Mathematik}, 39:119--137, 1982.

\bibitem{MoreGarbHill81}
J.~J. Mor\'{e}, B.~S. Garbow, and K.~E. Hillstrom.
\newblock Testing unconstrained optimization software.
\newblock {\em ACM Transactions on Mathematical Software}, 7(1):17--41, 1981.

\bibitem{PorcToin21}
M.~Porcelli and {Ph}.~L. Toint.
\newblock Exploiting problem structure in derivative-free optimization.
\newblock {\em ACM Transactions on Mathematical Software}, \textmd{(to
  appear)}, 2021.

\bibitem{Toin83a}
Ph.~L. Toint.
\newblock Test problems for partially separable optimization and results for
  the routine {PSPMIN}.
\newblock Technical Report 83/4, Department of Mathematics, FUNDP - University
  of Namur, Namur, Belgium, 1983.

\end{thebibliography}
\end{document}